\documentclass[a4paper,11pt]{article}
\usepackage[T1]{fontenc}
\usepackage[english]{babel} 
\usepackage[centertags,intlimits]{amsmath}
\usepackage{amsfonts}
\usepackage{amsthm}
\usepackage{amssymb}
\usepackage{mathrsfs}
\newtheorem{theorem}{Theorem}[section]
\newtheorem{proposition}[theorem]%
{Proposition}
\newtheorem{definition}{Definition}[section]
\newtheorem{lemma}[theorem]%
{Lemma}
\newtheorem{corollary}[theorem]%
{Corollary}

\newtheorem{example}{Example}[section]

\begin{document}

\title{Proximity and remoteness in triangle-free and $C_4$-free graphs in terms of order and minimum degree}

\author{
P.\ Dankelmann\thanks{Department of Mathematics and Applied Mathematics. University of Johannesburg. {\tt pdankelmann@uj.ac.za}.} \and E.\ Jonck\thanks{School of Mathematics. University of the Witwatersrand. {\tt Betsie.Jonck@wits.ac.za}.} \and S.\ Mafunda\thanks{Department of Mathematics and Applied Mathematics. University of Johannesburg. {\tt smafunda@uj.ac.za}.}}

\maketitle

\begin{abstract}
Let $G$ be a finite, connected graph. The average distance of a vertex $v$ of $G$ is 
the arithmetic mean of the distances from $v$ to all other vertices of $G$.  
The remoteness $\rho(G)$ and the proximity $\pi(G)$ of $G$ are the maximum and the 
minimum of the average distances of the vertices of $G$. 
In this paper, we present a sharp upper bound on the  remoteness of a triangle-free 
graph of given order and minimum degree, and a corresponding bound on the proximity, which
is sharp apart from an additive constant. We also present upper bounds on the remoteness
and proximity of $C_4$-free graphs of given order and minimum degree, and 
we demonstrate that these are close to being best possible. 
\end{abstract}

\section{Introduction}

Let $G$ be a finite, connected graph of order $n \geq 2$ with vertex set $V(G)$. 
The \textit{average distance} $\overline{\sigma}(v)$ of a vertex $v\in V(G)$ 
is defined as the arithmetic mean of the distances from $v$ to all other 
vertices of $G$, i.e.
$\overline{\sigma}(v,G)=\frac{1}{n-1} \sum_{w \in V(G)} d(v,w)$,
where $d(v,w)$ denotes the usual shortest path distance. 
The  {\em proximity} $\pi(G)$ is defined as $\min_{v\in V(G)} \overline{\sigma}(v)$, 
and the {\em remoteness} $\rho(G)$ of $G$ is defined as 
$\max_{v\in V(G)} \overline{\sigma}(v)$.

Bounds on proximity and remoteness in terms of order alone were given by 
Zelinka \cite{Zel1968} and later, independently, by Aouchiche and Hansen \cite{AouHan2010}, who
introduced the names proximity and remoteness.  

\begin{theorem}
{\rm (Zelinka \cite{Zel1968}, Aouchiche, Hansen \cite{AouHan2010})} \\
\label{theo:bounds-in-terms-of-order}
Let $G$ be a connected graph of order $n\geq 2$. Then
\[ \rho(G) \leq \frac{n}{2}, \]
with equality if and only if $G$ is a path and
\[ \pi(G) \leq  \left\{ \begin{array}{cc}
\frac{n+1}{4} & \textrm{if $n$ is odd,} \\
\frac{n+1}{4} + \frac{1}{4(n-1)} & \textrm{if $n$ is even.} 
\end{array} \right. \]
with equality if and only if $G$ is a path or a cycle. 
\end{theorem}

There are several results in the literature relating  
proximity or remoteness to other graph parameters, for example  
diameter \cite{Dan2018, AouHan2011}, radius \cite{WuZha2012, HuaCheDas2015} and 
average eccentricity \cite{MaWuZha2012}. For results related to proximity and
remoteness in trees see \cite{BarEntSze1997}. Remoteness in maximal 
planar graphs was considered in \cite{CzaDanOlsSze2019}.

The starting point for this paper is a strengthening of the bounds in 
Theorem \ref{theo:bounds-in-terms-of-order} that takes into account also 
the minimum degree. 

\begin{theorem} {\rm (Dankelmann \cite{Dan2015})}  \\
\label{theo:bounds-in-terms-of-min-degree}
Let $G$ be a connected graph of order $n$ and minimum degree $\delta$, 
where $\delta \geq 2$. Then
\[ \rho(G) \leq \frac{3n}{2(\delta+1)} +\frac{7}{2}, \]
\[ \pi(G) \leq \frac{3n}{4(\delta+1)} +3, \]
and both bounds are sharp apart from an additive constant. 
\end{theorem}

The goal of this paper is to show that the bounds in Theorem 
\ref{theo:bounds-in-terms-of-min-degree} can be strengthened significantly 
for triangle-free graphs, and also for graphs not containing a $4$-cycle.

\section{Terminology and Notation}

We use the following notation. 
Let $v$ be a vertex of $G$. Then the {\em neighbourhood} of $v$, denoted by
$N(v)$, is the set of all vertices adjacent to $v$. The {\em closed 
neighbourhood} $N[v]$ of $v$ is the set $N(v)\cup \{v\}$. 

The {\em total distance} (or {\em distance} for short) of $v$, $\sigma(v)$, 
is the sum of the distances from $v$ to all other vertices. Clearly
$\overline{\sigma}(v) = \frac{1}{n-1} \sigma(v)$.  
The {\em eccentricity} of $v$, denoted 
by ${\rm ecc}(v)$, is 
the distance from $v$ to a vertex farthest from $v$. 
The {\em radius} ${\rm rad}(G)$ of $G$ is the smallest of all eccentricities of the vertices
of $G$. A vertex whose eccentricity equals ${\rm rad}(G)$ is called a {\em centre vertex}
of $G$. 
If $i\in\mathbb{Z}$, then 
$N_i(v)$ is the set of all vertices at distance $i$ from $v$, and $n_i$ its cardinality.
By $N_{\leq i}(v)$ and $N_{\geq i}(v)$ we mean the set of vertices at distance
at most $i$ and at least $i$, respectively, from $v$. 
Clearly we have $n_i \geq 1$ if and only if $0\leq i\leq {\rm ecc}(v)$. 
By the {\em distance degree} $X(v)$ of $v$ 
we mean the sequence $(n_0,n_1,\ldots,n_d)$, where $d$ is the 
eccentricity of $v$. 

Here, and for any finite sequence $(n_0,\ldots,n_d)$, we use 
the convention that $n_i=0$ for all $i\in \mathbb{Z}$ with $i<0$ or $i>d$. 
With this convention, we define for a finite sequence $X=(n_0,n_1,\ldots,n_d)$,  
\[g(X)=\sum\limits_{i=0}^{\infty}in_i.\]
So if $X(v)=(n_0,n_1,\ldots,n_d)$, then 
clearly  $\sigma(v) = g(X(v))$. Note that we will generally apply the function $g$ 
only to sequences that have a finite number 
of non-zero entries, so questions of convergence won't arise.
If $X$ and $X^{*}$ are two sequences of nonnegative integers, then we say that 
$X^{*}$ {\em beats}  $X$ if $g(X^*)>g(X)$. 

We often modify a given sequence $X=(n_0, n_1,\ldots,n_d)$. 
If $a\in \mathbb{N}\cup \{0\}$ and $X^*$ is the sequence obtained from $X$ by adding
$a$ to $n_i$, i.e. if 
$X^*=(n_0,\ldots,n_{i-1}, n_i+a, n_{i+1}\ldots,n_d)$, then we say that $X^*$ is
obtained from $X$ by applying $n_i \leftarrow +a$. Similarly we define the
sequence $X^*$ obtained from $X$ by applying $n_i \leftarrow -a$.
We usually write $X^*$ as $(n_0', n_1',\ldots,n_{d'}')$. 
For a given sequence $(n_0,\ldots,n_d)$ we also define 
$S_3(i)=n_{i-1}+n_i+n_{i+1}$, $S_4(i)=n_{i-1}+n_i + n_{i+1} + n_{i+2}$
and $S_5(i)=n_{i-2}+n_{i-1}+n_i+n_{i+1}+n_{i+2}$ for $i\in \mathbb{Z}$. 
We write $S_3'(i)$, $S_4'(i)$ and $S_5'(i)$ for the corresponding terms 
of a modified sequence $X^*$.

By $K_n$, $\overline{K_n}$, $C_n$ and $P_n$ we mean the complete graph, the edgeless graph, the 
cycle, and the path on $n$ vertices. By a {\em triangle} we mean the graph $K_3$.
If $F$ is a graph, then we say that $G$ is $F$-{\em free} if $G$ does not contain $F$
as a (not necessarily induced) subgraph.

\section{Remoteness in triangle-free graphs}

In this section we give an upper bound on the remoteness of a 
triangle-free graph in terms of order and minimum degree. Unlike the bound
in Theorem \ref{theo:bounds-in-terms-of-min-degree}, our bound is sharp.   
We prove our bound by first demonstrating some properties of the distance 
degree $X(v)$ of an arbitrary vertex $v$ of a triangle-free graph, and 
then determining a sequence that maximises the function $g$ among all 
sequences with these properties.

\begin{proposition} \label{prop:distance-seq-triangle-free}
Let $u$ be a vertex of a connected, triangle-free graph of order $n$ and minimum degree $\delta$. Let $X(u)=(n_0,n_1,\ldots,n_d)$. Then the following hold: \\
(A1) $n_0=1$, \\
(A2) $\sum\limits_{i=0}^\infty n_i=n$, \\
(A3) if $i\geq 1$ and $n_i\geq 1$, then $n_1, n_2, \cdots, n_{i-1}\geq 1$, \\
(A4) if $i\geq 0$ and $n_i, n_{i+1}\geq 1$, then 
$S_4(i) \geq 2\delta$, \\
(A5) if $i\geq 0$ and $n_i>0$, then $S_3(i) \geq \delta +1$, \\
(A6)  if $n_{d-1} \leq \delta-1$, then $n_{d-1} + n_d \geq 2\delta$. 
\end{proposition}

{\bf Proof:} Let $N_i=N_i(v)$ and $n_i=n_i(v)$. 
Clearly, (A1), (A2) and (A3) hold. \\
(A4): Assume that $n_i, n_{i+1}\geq 1$.
Let $v_{i+1}$ be a vertex in $N_{i+1}$. Then $v_{i+1}$ is adjacent to a vertex $v_i\in N_i$. 
Since, $G$ is triangle-free, $N_G(v_i)$ and $N_G(v_{i+1})$ are disjoint subsets of the vertex 
set, we have 
$|N(v_i) \cup N(v_{i+1})| = \deg_G(v_i) + deg_G(v_{i+1})\geq 2\delta$. 
Since
$N_G(v_i)\cup N_G(v_{i+1})\subseteq N_{i-1}\cup N_{i}\cup N_{i+1}\cup N_{i+2}$, 
we have  
\[ S_4(i) = n_{i-1} + n_{i} + n_{i+1} + n_{i+2} \geq |N_G(v_i)|+|N_G(v_{i+1})| \geq 2\delta, \]
and (A4) follows. \\
(A5): Assume that $n_i>0$. Let $v_i \in N_i$. Since 
$|N_G[v_i]| ={\rm deg}_G(v_i) + 1 \geq \delta+1$ 
and  $N_G[v_i]\subset N_{i-1}\cup N_i\cup N_{i+1}$, we conclude 
that $n_{i-1}+n_i+n_{i+1}\geq \delta +1$, and (A5) follows. \\
(A6): Assume that $n_{d-1} \leq \delta-1$. Let $v_d \in N_d$. Since 
$N(v_d) \subseteq N_{d-1} \cup N_d$, and since $n_{d-1}<\delta$, vertex $v_d$ has a 
neighbour $v_d'$ in $N_d$. Since $G$ is triangle-free, $N(v_d)$ and $N(v'_d)$ are disjoint
and so $|N(v_d) \cup N(v_d')| = |N(v_d)| + |N(v_d')| \geq 2\delta$. 
Since $N(v_{d}) \cup N(v_{d}') \subseteq N_{d-1} \cup N_{d}$, we conclude 
that $n_{d-1} + n_d \geq 2\delta$, which is (A6).
\hfill $\Box$

\noindent 
For the following definition we need to introduce some notation. 
Given $p\in \mathbb{N}$ and an infinite sequence $A=(a_1, a_2,\ldots)$ of positive
integers, we define $\ell(A,p)$ to be the smallest value $k$ with 
$\sum_{i=1}^k a_i > p$. 

\begin{definition}\label{def3.2.1}
Given $n,\delta \in \mathbb{N}$ with $\delta\geq 3$. Let 
$A=(1,1,\delta-1, \delta-1, 1, 1, \delta-1, \delta-1,\ldots)$ be the 
infinite sequence repeating the $(1,1,\delta-1,\delta-1)$-pattern indefinitely. 
Define the finite sequence $X_{n,\delta}$ by 
\[ X_{n,\delta} = (1, \delta, \delta-1, a_1, a_2,\ldots,a_{\ell(A,n-4\delta)}, \delta, r_{n,\delta}), \]
where $r_{n,\delta} = n-3\delta - \sum_{i=1}^{\ell(A,n-4\delta)} a_i$.    
\end{definition}

We now show that, for given $n$ and $\delta$, a sequence that satisfies (A1)-(A6), 
and that is not beaten
by any other such sequence, necessarily equals $X_{n,\delta}$.

\begin{lemma}\label{lem3.2.4}
Let $n,\delta \in \mathbb{N}$ with $n\geq 6\delta$ and $\delta\geq 2$.
If $X=(n_0, n_1,\ldots,n_d)$ is a sequence satisfying (A1)-(A6) that is not beaten
by any other such sequence, then $X=X_{n,\delta}$. 
\end{lemma}

{\rm Proof:} 
Let $X$ be a sequence satisfying (A1)-(A6). It suffices to prove the lemma for 
sequences which is not beaten by any other sequence 
satisfying (A1)-(A6), so let $X$ be such a sequence. 
Note that only $n$ and $\delta$ are fixed, but $d$ is not. 
In a sequence of claims we prove properties of the sequence $X$ 
by showing that otherwise we can modify $X$ to obtain a sequence 
$X^*= (n_0',n_1',\ldots,n_{d'}')$ which satisfies (A1)-(A6) 
but beats $X$, thus obtaining a contradiction.  \\[1mm]
{\sc Claim 1:} If $2 \leq i \leq d-2$, then $n_{i}\leq \delta-1$. \\
Suppose to the contrary that for some $j\in \{2,\ldots, d-2\}$, we 
have $n_j\geq\delta$. Let $X^*$ be the sequence obtained from $X$
by applying $n_j\leftarrow-1$  and   $n_{j+1}\leftarrow+1$. 
Then clearly $X^*$ satisfies conditions (A1), (A2) and (A3). 
We now show that (A4) holds for $i=j-2$. Clearly, $n_{j-3}', n_{j-2}', n_{j-1}', n_{j}'$ is the 
only quadruple in $X^*$ where possibly condition (A4) may fail since only the values 
$n_j$ and $n_{j+1}$ have changed. Now 
$n_{j-3} + n_{j-2} + n_{j-1}\geq \delta +1$ by (A5), hence 
$n_{j-3}' + n_{j-2}' + n_{j-1}' + n_{j}' 
  = (n_{j-3} + n_{j-2} + n_{j-1}) + (n_{j}-1) \geq (\delta+1)+(\delta-1)=2\delta$, 
so (A4) holds for $i=j-2$, and thus for all $i\in\{3, \dots, d-2\}$. We now show 
that (A5) holds. Since only the values $n_j$ and $n_{j+1}$ have changed, 
condition (A5) holds for all $i$ except possibly $i=j-1$. But $n_{j-2}, n_{j-1}\geq 1$, 
so $n_{j-2}' + n_{j-1}' + n_{j}' \geq n_{j-2} + n_{j-1} + (n_{j}-1) \geq \delta+1$. 
Hence condition (A5) holds for $i=j-1$ and thus for all $i$. Condition (A6) 
holds since $d'=d$, $n'_{d-1}\geq n_{d-1}$ and $n'_{d} = n_{d}$. 
Since $g(X^*)=g(X)+1$, $X^*$ beats $X$, a contradiction to the choice of $g$. 
This proves Claim 1. \\[1mm]
{\sc Claim 2:} $(n_0, n_1, n_2, n_3, n_4, n_5)=(1, \delta, \delta -1, 1, 1, \delta -1)$. \\
We determine the values of $n_0$ to $n_5$ 
in five subclaims. \\[1mm]
{\sc Claim 2.1:} $n_0=1$ and $n_1=\delta$. \\
Clearly, $n_0=1$ by (A1). Condition (A5) for $i=0$ yields $n_1 \geq \delta$. 
Suppose to the contrary that $n_1>\delta$. Let $X^*$ be the sequence obtained 
from $X$ by applying $n_1\leftarrow-1$ and $n_2\leftarrow+1$.  
It is easy to verify that $X^{*}$ satisfies (A1)-(A6). Moreover, 
$g(X^*) = g(X)+ 1 > g(X)$, so $X^{*}$ beats $X$. This contradiction proves 
Claim 2.1. \\[1mm]
{\sc Claim 2.2:} $n_2=\delta-1$. \\
It follows from (A4) for $i=0$ that $n_{2}\geq\delta-1$. On the other
hand we have $n_2 \leq \delta-1$ by Claim 1.  Claim 2.2 follows. \\[1mm]
{\sc Claim 2.3:}  $n_3=1$.\\
By condition (A3) we have  $n_3\geq 1$. Suppose to the contrary that 
$n_3>1$. Let $X^*$ be the sequence obtained from $X$ by applying 
$n_3\leftarrow-1$ and $n_4\leftarrow+1$.  
Then clearly $X^*$ satisfies (A1), (A2) and (A3). We now show that (A4) holds. 
Clearly (A4) holds for all $i$ except possibly $i=1$ since only the values $n_3$ and 
$n_4$ have changed. Consider  $i=1$, we have 
$n_0'+n_1'+n_2'+n_3'=1 + \delta+(\delta -1)+(n_3-1)> 2\delta$, 
which implies that (A4) holds for $i=1$ and thus for all $i$. Next we show that 
(A5) holds. To do so, consider $i=2$ since only the values $n_3$ and $n_4$ have changed. 
Since $n_1'+n_2'+n_3'=\delta+(\delta -1)+(n_3-1)\geq 2\delta >\delta +1$, 
condition (A5) holds for all $i$. It remains only to show that (A6) is satisfied. 
Since  $n_{d'-1}'\geq n_{d-1}$ and $n_{d'}'=n_d$, (A6) holds. 
Since $g(X^*)=g(X)+1$, $X^*$ beats $X$. This contradiction to the 
choice of $X$ proves Claim 2.3. \\[1mm]
{\sc Claim 2.4:} $n_4=1$.\\
By condition (A5) for $i=3$, we have that $n_4\geq 1$. Suppose to the 
contrary that $n_4>1$. Now, if $X^*$ is a sequence obtained from $X$ by applying 
$n_4\leftarrow-1$ and $n_5\leftarrow+1$.  
Arguments similar to those in Claim 2.3 show that $X^*$ satisfies (A1)-(A6).  
Since $g(X^*)=g(X)+1$, $X^*$ beats $X$, a contradiction. Hence $n_4=1$. \\[1mm] 
{\sc Claim 2.5:} $n_5=\delta -1$.\\
By condition (A5) for $i=4$, we have $n_5\geq \delta -1$. Also, 
from Claim 1 we have that $n_5\leq \delta-1$. Hence $n_5=\delta -1$, which
is Claim 2.5. \\[1mm]
{\sc Claim 3:} Let $3 \leq i \leq d-4$, then $S_4(i)=2\delta$. \\
From (A4) we have that  
$S_4(i) \geq 2\delta$ 
for all $i\in \{3,4,\ldots,d-4\}$. 
Suppose to the contrary that there exists $i\in \{3,4,\ldots,d-4\}$ 
such that $S_4(i) > 2\delta$.
Choose a smallest such $i$. It follows from Claim 2 that 
$S_4(3)=n_2+n_3+n_4+n_5=2\delta$, so $i>3$.
Let $X^{*}$ be the sequence obtained from $X$ by applying
$n_{i+2}\leftarrow -1$ and  $n_{i+3}\leftarrow +1$. 
Then, clearly $X^{*}$ satisfies conditions (A1) and (A2). 
Since $n_{i-1}, n_i, n_{i+1}, n_{i+2}$ is the first quadruple in $X$ for which 
$n_{i-1}+n_i+n_{i+1}+n_{i+2}> 2\delta$, we have that $n_{i+2}>1$ since by 
the minimality of $i$ and by $i>3$ we have $n_{i-2} + n_{i-1} + n_i + n_{i+1} = 2\delta$, 
but $n_{i-1} + n_i + n_{i+1} + n_{i+2} > 2\delta$, implying that 
$n_{i+2} > n_{i-2} \geq 1$. Hence $X^*$ satisfies condition (A3). 
Also, $n_{i-1} + n_i + n_{i+1} + n_{i+2} > 2\delta$ implies 
$n_{i-1}+n_i+n_{i+1}+(n_{i+2}-1)\geq 2\delta$. Since $n_{i-1}', n_i', n_{i+1}', n_{i+2}'$ 
is the only quadruple whose sum in $X^*$ is less than that in $X$,
$X^*$ satisfies (A4).\\
Now we show that (A5) holds. Suppose to the contrary that there exists a 
$j$ for which $n_{j-1}'+n_{j}'+ n_{j+1}' \leq \delta$. Then $j=i+1$ since 
$n_i, n_{i+1}, n_{i+2}$ is the only triple in $X^*$ whose sum is less than in $X$. Since 
$n_{i}+n_{i+1}+n_{i+2} \geq \delta+1$, we conclude that 
$n_i + n_{i+1} + n_{i+2}=\delta+1$.  Since $n_{i-1}+n_i+n_{i+1}+n_{i+2} \geq 2\delta+1$, we have $n_{i-1} \geq \delta$, which contradicts ($P1$), hence sequence $X^*$ 
satisfies (A5). The sequence $X^*$ satisfies also (A6) since $i\leq d-4$ implies that 
$d'=d$,  $n_{d'-1}'\geq n_{d-1}$ and $n_{d'}'= n_d$. 
Moreover, $g(X^*) = g(X)+1 > g(X)$,   
so $X^{*}$ beats $X$. This contradiction to the assumption that $X$ is in 
not beaten completes the proof of Claim 3.  \\[1mm]
{\sc Claim 4:} $n_{d-1}=\delta$ and $n_d \leq \delta-1$. \\
We consider three cases, depending on the value of $n_d$. \\[1mm] 
{\sc Case 1}:  $n_d \geq  \delta+1$. \\
Let $X^*$ be the sequence obtained from $X$ by applying $n_d \leftarrow -1$
and $n_{d+1} \leftarrow +1$, so $d'=d+1$. Then it is clear that $X^*$ satisfies 
(A1), (A2) and (A3). We now show that (A4) holds. 
Clearly, (A4) holds for all $i$ except possibly $i=d-2$ or $i=d$ since only
the values $n_d$ and $n_{d+1}$ have changed. For $i=d-2$ we have 
$n_{d-3}' + n_{d-2}' + n_{d-1}' + n_d' = (n_{d-3} + n_{d-2} + n_{d-1}) + (n_d-1) 
       \geq (\delta+1) + \delta > 2\delta$ by (A5),
and for $i=d$ we have 
$n_{d-1}' + n_d' + n_{d+1}' + n_{d+2}' = n_{d-1}+n_d \geq 2\delta$ 
since either $n_{d-1} \geq \delta$, which implies 
$n_{d-1}+n_d \geq \delta + \delta+1 >2\delta$,
or $n_{d-1} \leq \delta-1$, which implies $n_{d-1}+n_d \geq 2\delta$ by (A6).
Hence (A4) holds. 
We now show that (A5) holds.  
Since only the values $n_d$ and $n_{d+1}$ have changed, condition (A5) holds for all
$i$ except possibly for $i=d+1$ or $i=d-1$. Since 
$n_d' + n_{d+1}' + n_{d+2}' = n_d \geq \delta+1$, condition (A5) holds for $i=d+1$, 
and since $n_{d-2}'+ n_{d-1}' + n_d' \geq 1+1+(n_d-1) \geq \delta+1$, 
condition (A5) holds for $i=d-1$. Hence (A5) holds for all $i$. 
Condition (A6) holds vacuously since $d'=d+1$ and $n_{d'-1}'=n_d-1 \geq  \delta$.
Since clearly $g(X^*) = g(X)+1$, $X^*$ beats $X$.  This contradiction to the choice 
of $X$ proves that Case 1 cannot occur. \\[1mm]
{\sc Case 2:} $n_d =  \delta$. \\
Then (A6) implies that $n_{d-1} \geq \delta$. 
Let $X^*$ be the new sequence obtained from $X$ by applying $n_{d-1} \leftarrow -1$
and $n_{d+1} \leftarrow +1$. Clearly, $X^*$ satisfies (A1), (A2) and (A3).  
We now show that (A4) holds. 
Clearly, (A4) holds for all $i$ except possibly $i=d-3$, $i=d-2$ or $i=d$ since only
the values $n_{d-1}$ and $n_{d+1}$ have changed. Now 
$n_{d-4}' + n_{d-3}' + n_{d-2}' + n_{d-1}' = (n_{d-4} + n_{d-3} + n_{d-2})) + (n_{d-1}-1) 
       \geq (\delta+1) + (\delta-1) = 2\delta$,
and 
$n_{d-3}' + n_{d-2}' + n_{d-1}' + n_d' = n_{d-3} + n_{d-2} + (n_{d-1}-1) + n_d 
       \geq 1 + 1 + (\delta-1) + \delta \geq 2\delta$,
so (A4) holds for $i\in \{d-3,d-2\}$.        
Since 
$n_{d-1}' + n_d' + n_{d+1}' + n_{d+2}' = n_{d-1}+n_d \geq 2\delta$,  
(A4) holds also for $i=d$, and thus for all $i$. 
We now show that (A5) holds.  
Since only the values $n_{d-1}$ and $n_{d+1}$ have changed, condition (A5) holds for all
$i$ except possibly for $i=d-2$, $i=d-1$, or $i=d+1$. 
Since $n_{d-1}' \geq \delta-1$ and $n_{d-3}',n_{d-2}', n_d'\geq 1$, 
condition (A5) holds for $i=d-2$ and for $i=d-1$. Since $n_{d}'+n_{d+1}'=\delta+1$,
it also holds for $i=d+1$. 
Condition (A6) holds vacuously since $n_{d'-1}=n_d=\delta$. 
Since clearly $g(X^*) = g(X)+1$, $X^*$ beats $X$. This contradiction to the choice of
$X$ shows that Case 2 cannot occur. \\[1mm] 
{\sc Case 3}: $n_d \leq \delta-1$. \\
In this case we only have to show that $n_{d-1}=\delta$. 
Suppose not. It follows from (A6) that $n_{d-1} \geq \delta$,
so $n_{d-1} \geq \delta+1$. 
Let $X^*$ be the new sequence obtained from $X$ by applying $n_{d-1} \leftarrow -1$
and $n_{d} \leftarrow +1$. Then as above, we show that $X^*$
beats $X$, a contradiction. Hence $n_{d-1}=\delta$, which completes 
the proof of Claim 4. 

We are now in a position to complete the proof of Lemma \label{lem3.2.4} by
showing that $X=X_{n,\delta}$. Let  $X_{n,\delta}=(x_0,x_1,\ldots,x_{\ell})$. 

Claim 2 implies that $n_i=x_i$ for $0 \leq i \leq 5$. Since Claim 3 holds for
$X$ and for $X_{n,\delta}$ (with $d$ replaced by $\ell$), we have 
$n_i=x_i$ for all $i\in \{6,7,\ldots,\min(d-2, \ell-2) \}$. 
We now consider three cases. \\[1mm]
{\sc Case 1:} $\ell=d$. \\
By Claim 4 and $d=\ell$ we have $n_{d-1}=x_{d-1}$, 
and so, by $\sum_{i=0}^d n_i = \sum_{i=0}^d x_i$ we also have
$n_d=x_d$, which implies $X=X_{n,d}$, as desired. \\[1mm]
{\sc Case 2:} $\ell<d$. \\
Then $n_i=x_i$ for $0 \leq i \leq \ell-2$. Hence 
$x_{\ell-1} + x_{\ell}= n_{\ell-1} + n_{\ell} + \cdots + n_d \geq 2\delta$
by (A4). On the other hand, we have $x_{\ell-1}=\delta$ and 
$x_{\ell}=r_{n,\delta}=n-3\delta - \sum_{i=1}^{\ell(A,n-4\delta)} a_i 
         < n-3\delta- (n-4\delta) = \delta$, 
so $x_{\ell-1} + x_{\ell} < 2\delta$. 
This contradiction proves that Case 2 cannot occur. \\[1mm]
{\sc Case 3:} $\ell>d$. \\
By Claim 4 we have $n_{d-1} + n_d \leq 2\delta-1$. Since 
$x_{d-1} + x_d  \cdots + x_{\ell} = n_{d-1} + n_d$ and $x_{\ell-1}=\delta$ we have 
$x_{\ell-2} + x_{\ell} \leq \delta-1$. By the definitions of $X_{n,\delta}$ 
and $A$ we have $x_{\ell-2} = a_{\ell(A, n-4\delta)} \in \{1, \delta-1\}$, so 
$x_{\ell-2} = a_{\ell(A, n-4\delta)}=1$. Hence it follows from the definition of 
$\ell(A, n-4\delta)$ that 
$x_3 + x_4 + \cdots, x_{\ell-2} = \sum_{i=1}^{\ell(A,n-4\delta)} a_i= n-4\delta+1$. 
Hence $x_0 + x_1 + \cdots, x_{\ell-2} = n-2\delta+1$, and, since $x_{\ell-1}= \delta$,  
also $x_{\ell} = \delta-1$, a contradiction to $x_{\ell-2} + x_{\ell} \leq \delta-1$. This 
contradiction, which proves that Case 3 cannot occur, completes the proof
of the lemma. 
\hfill $\Box$

For a finite sequence $X=(x_0,x_1,\ldots,x_d)$ of positive integers we define
the graph $G(X)$ by 
\[ G(X) = \overline{K_{x_0}} +  
    \overline{K_{x_1}} +  \ldots +
      \overline{K_{x_d}}. \]

\begin{theorem} \label{theo:bound-on-remoteness-triangle-free}
Let $G$ be a connected, triangle free graph of order $n$ and minimum degree 
$\delta$, where $\delta \geq 3$ and $n \geq 6\delta$. Then 
\[ \rho(G) \leq \rho(G(X_{n,\delta})). \]
\end{theorem}  

{\bf Proof:}
Let $n$ and $\delta$ be fixed. Assume $G$ is a connected triangle-free graph of order $n$ and minimum degree $\delta$. Let $u$ be a vertex of maximum distance in $G$, and let $X(u)=(n_0, n_1,\ldots,n_d)$  
be its distance degree. Then
\[ \rho(G) = \overline{\sigma}(u)= \frac{1}{n-1}g(X(u)). \]
Since $X(u)$ satisfies (A1)-(A6), and since among all such sequences 
$X_{n,\delta}$ maximises $g$ by Lemma \ref{lem3.2.4}, we have  
\[ \frac{1}{n-1} g(X(u)) \leq   \frac{1}{n-1} g(X_{n,\delta}). \]
If $v$ is the vertex of $G(X_{n,\delta})$ contained in $K_{x_0}$, then 
the distance degree of $v$ is $X_{n,\delta}$, so 
\[ \frac{1}{n-1} g(X_{n,\delta}) 
  = \overline{\sigma}(v,G(X_{n,\delta})) 
  \leq \rho(G(X_{n,\delta})). \] 
Combining the last three (in)equalities yields the theorem. 
\hfill $\Box$

It is easy to verify that $G(X_{n,\delta})$ is a bipartite 
(and thus triangle-free) graph
of order $n$ and minimum degree $\delta$. Hence the bound in Theorem
\ref{theo:bound-on-remoteness-triangle-free} is sharp. 
Evaluating the remoteness of $G(X_{n,\delta})$ yields the 
following corollary.

\begin{corollary}\label{thm3.2.5}
Let $n,\delta\in\mathbb{N}$, with $\delta\geq 3$. If $G$ is a connected, triangle-free graph of order $n$ and minimum degree $\delta$, then
\[\rho(G)\leq 2\Big\lceil\frac{n-3\delta}{2\delta}\Big\rceil +2-\frac{\delta}{n-1}\] 
and this bound is sharp.
\end{corollary}

\section{Proximity in triangle-free graphs}
\label{section:proximity-in-triangle-free}

In this section we present an upper bound on the proximity of triangle-free 
graphs of given order and minimum degree which improves on the bound in
Theorem \ref{theo:bounds-in-terms-of-min-degree} and is sharp apart
from an additive constant. 

Throughout this section let $G$ be a connected, triangle-free graph of order $n$, 
minimum degree $\delta$, and radius $r$. Furthermore let $v_0$ be a fixed central vertex 
of $G$, let $v_r$ be a fixed vertex at distance $r$ from $v_0$, and let $T$ be a 
spanning tree of $G$ that preserves the distances from $v_0$.
For vertices $a$ and $b$ of $T$ we denote the unique $(a,b)$ path in $T$ by $T(a,b)$. 
By $N_i$ and $n_i$ we mean $N_i(v_0)$ and $n_i(v_0)$, respectively.

The following two propositions show that the 
the properties of the distance degree of an arbitrary vertex 
in Proposition \ref{prop:distance-seq-triangle-free} can be strengthened for centre vertices.

\begin{proposition} \label{prop:ni>=2-for-centre-vertex}
Let $v_0$ be a centre vertex of a connected graph $G$ of radius $r$ on at least
three vertices. Then 
$n_i \geq 2$ for all $i \in \{1,2,\ldots,r-1\}$. 
\end{proposition}

{\bf Proof:}
Clearly, a vertex of degree $1$ cannot be a centre vertex since its neighbour 
has smaller eccentricity. Hence we have $n_1 \geq 2$. 

We now prove that $n_i \geq 2$ for all $2 \leq i \leq r-1$. 
Suppose to the contrary that for some $i\in\{2,\ldots, r-1\}$ we have $n_i=1$. 
Consider a vertex $v_1$ of $N_1$ on a shortest $(v_0,v_i)$-path. 
We obtain a contradiction by showing that ${\rm ecc}(v_1)<{\rm ecc}(v_0)$.
Let $w_j$ be a vertex in $N_j$ for $0\leq j\leq i-1$. Then 
\[ d(v_1,w_j) \leq d(v_1, v_0) + d(v_0, w_j)
               \leq 1 + i
               \leq r-1.   \]
Now consider a vertex $w_k \in N_k$ for $i\leq k\leq r$. 
Since $v_1$ is on a shortest $v_0-v_i$ path, and since every $v_0-w_k$ path 
goes through $v_i$, we have 
\[
d(v_1, w_k) = d(v_1,v_i) + d(v_i, w_k)
    = (i-1) + (k-i)
    = k-1
    \leq r-1.   \] 
If follows that ${\rm ecc}(v_1) < r= {\rm rad}(G)$, a contradiction. 
Hence $n_i\geq 2$ for any $2\leq i\leq r-1$. \hfill $\Box$

Following \cite{ErdPacPolTuz1989}, we say that two vertices $v$ and $v'$ of $G$ are {\em  related} if 
there exist vertices  
$x \in V(T(v_0,v)) \cap N_{\geq 9}(v_0)$ and $y \in V(T(v_0,v')) \cap N_{\geq 9}(v_0)$ 
such that $d_G(x,y) \leq 4$.

\begin{lemma}  {\rm (Erd\"{o}s, Pach, Pollack, Tuza \cite{ErdPacPolTuz1989})} 
\label{la:related}   \\
Let $G$ be a connected graph of radius $r$. Let $v_{0}$ be a centre vertex of $G$ and $T$ a spanning tree of $G$ that preserves the distances from $v_0$. If $v_r\in N_r$, then there exists a vertex in $N_{\geq r-9}$ which is not related to $v_r$.  
\end{lemma}

\begin{proposition}\label{prop3.2.6}
Let $G$ be a connected triangle-free graph of order $n$, radius $r$ 
and minimum degree $\delta \geq 3$. If $v_0$ is a centre vertex of $G$,
then the distance degree $X(v_0)=(n_0,n_1,\ldots,n_r)$ satisfies 
the following. 
\begin{enumerate}
\item[(B1)] $n_0=1$,
\item[(B2)] $\sum\limits_{i=0}^\infty n_i=n$,
\item[(B3)] If $i>0$ and $n_i>0$, then $n_1, n_2,\ldots, n_{i-1} \geq 2$,
\item[(B4)] $S_4(i) \geq 2\delta$ for all $i\in\{1,\ldots,r-1\}$ with $i\equiv 1\pmod{4}$,
\item[(B5)] $S_4(i) \geq 4\delta$ for all $i\in\{9,\ldots,r-10\}$ with $i\equiv 1 \pmod{4}$
\end{enumerate}
\end{proposition}

{\bf Proof:} 
(B1), (B2) and (B4) are the properties (A1), (A2) and (A4), respectively, in 
Proposition \ref{prop:distance-seq-triangle-free}.  
(B3) follows from Proposition \ref{prop:ni>=2-for-centre-vertex}. \\[1mm]
(B5): Let $v_r\in N_r$.   
By Lemma \ref{la:related} there is a vertex $v_k\in N_k$, for some $k\geq r-9$, that is not related to $v_r$. 
Define 
\[R=V\left(T(v_{0},v_r)\right)\cap N_{\geq 9}\qquad\mbox{and}\qquad K=V\left(T(v_{0},v_k)\right)\cap N_{\geq 9}.\] 
For $i\in \{0,1,\ldots,d\}$ let  
\[N_i^{\prime}=\{v_i\in N_i\;|\; d_G(v_i,R)\leq 2\}\qquad\mbox{and}\qquad N_i^{\prime\prime}=\{v_i\in N_i\;|\; d_G(v_i,K)\leq 2\}.\] 
Since $v_k$ is not related to $v_r$, we have $d_G(x,y) \geq 5$ for all
$x\in R$ and $y\in K$. Hence $N_i^{\prime}\cap N_i^{\prime\prime}=\emptyset$ 
for $i=7,8,\ldots,\min(k+2,r)$.
Letting $n_i^{\prime}=|N_i^{\prime}|$ and 
$n_i^{\prime\prime}=|N_i^{\prime\prime}|$, we have 
$n_i\geq n_i^{\prime}+n_i^{\prime\prime}$ for all $i$ with 
$7\leq i\leq \min(k+2,r)$.
Let $9 \leq i \leq r-10$. Let $v_i'$ and $v_{i+1}'$ be the unique vertices of $R$ that are in $N_i$ and $N_{i+1}$ respectively. Then $N_G(v_i')$ and $N_G(v_{i+1}')$ are contained in $N_{i-1}' \cup N_i'\cup N_{i+1}'\cup N_{i+2}'$. Since $G$ is triangle-free these neighbourhoods are disjoint. Hence, since ${\rm deg}_G(v_i') \geq \delta$ and ${\rm deg}_G(v_{i+1}') \geq \delta$, we have 
\[ 2\delta \leq 
    |N_G(v_i') \cup N_G(v_{i+1}')| 
   \leq n_{i-1}'+n_i'+ n_{i+1}'+ n_{i+2}'. \]
Similarly, if $v_i''$ and $v_{i+1}''$ are the unique vertices of $K$ that are in $N_i$ and $N_{i+1}$, respectively, we obtain
\[ 2\delta \leq 
    |N_G(v_i'') \cup N_G(v_{i+1}'')| 
   \leq n_{i-1}''+n_i''+ n_{i+1}''+ n_{i+2}''. \]
Hence, since $n_j \geq n_j'+ n_j''$ for all $j\in \{0,1,\ldots,d\}$, 
\[ S_4(i) = n_{i-1}+n_i+n_{i+1}+n_{i+2}\geq 4\delta,\quad\mbox{for}\quad 9\leq i\leq r-10.\]
Since $n_{i-1} + n_i +n_{i+1} + n_{i+2} \geq 4\delta$ holds for all $i$ with $9 \leq i \leq r-10$, it holds in particular for all such $i$ with $i\equiv 1 \pmod{4}$. Note that we will find it convenient later to consider this property only for values of $i$ that are congruent to $1 \pmod{4}$. 
\hfill $\Box$

\begin{definition}\label{def3.2.2}
(a) Given $n, \delta \in \mathbb{N}$ with $\delta \geq 4$ and $n > 15\delta+3$.  
Let $p=\lceil\frac{n-15\delta-3}{4\delta}\rceil$ and $n_r=n-(4p+8)\delta-4$. 
We define 
$Y_{n,\delta}$ to be the finite sequence 
\[  (1, 2, 2, 2\delta-5, 2,2,2,2\delta -6, [2,2,2,4\delta -6]^{p}, 
                2,2,2,2\delta -6, 2,2,2,2\delta -6, 2, 2, n_r),\]
where $[2,2,2,4\delta-6]^{p}$ stands for the $p$-fold repetition of the 
quadruple $(2,2,2,4\delta-6)$.  \\
(b) Given $n\in \mathbb{N}$ with $n > 51$. Let $p=\lceil\frac{n-54}{12}\rceil$ 
and $n_r=n-12p-35$. We define $Y_{n,3}$ to be the finite sequence  
\[ (1, 2, 2, 2, 2,2,2,2, [2,2,2,6]^{p}, 2,2,2,2, 2,2,2,2, 2, 2, n_r),\]
where $[2,2,2,6]^{p}$ stands for the $p$-fold repetition of the 
quadruple $(2,2,2,6)$.          
\end{definition}

A long but straightforward calculation shows that  for $\delta\geq 4$
\begin{eqnarray*} 
g(Y_{n,\delta}) &=& g((1,2,2,2\delta-5)) + g(([0]^4, 2,2,2,2\delta-6)) \\
   & & + g(([0]^8,[2,2,2,4\delta-6]^p, 2,2,2,2\delta-6,2,2,2,2\delta-6,2,2,n_r)) \\
   &=& [20\delta -21] + [-12p+ 8\delta p^2 +36\delta p]\\ 
     & & + [16\delta p + 16p + 4pn_r +42 + 52\delta + 18n_r] \\
   &=& 8\delta p^2 + 52\delta p + 4p + 4pn_r  + 72\delta + 21 + 18n_r]
 \end{eqnarray*}

In the following we employ a similar proof strategy as for the bound on remoteness 
in triangle-free graphs.  We prove that $Y_{n,\delta}$ is not beaten by any 
sequence satisfying (B1)-(B5), and so we show that  
$g(Y_{n,\delta})$ is an upper bound for the total distance of vertex $v_0$.
Some details, however, are more involved. For example, it turns out that for some values 
of $n$ and $\delta$ the sequence $Y_{n,\delta}$ is not the only sequence maximising $g$.

\begin{lemma}\label{la:maximal-sequence-for-proximity-in-triangle-free}
Let $n, \delta \in \mathbb{N}$ with $n > 15\delta+3$ and $\delta \geq 4$ be given. 
If $X$ is a sequence that satisfies (B1)-(B5), then 
\[ g(X) \leq g(Y_{n,\delta}). \]
\end{lemma}

{\bf Proof:}
It suffices to prove the lemma for sequences satisfying (B1)-(B5) that are not beaten
by any other sequence satisfying (B1)-(B5). Let $X=(n_0,\ldots,n_r)$ be such a sequence. 
Note that $n$ and $\delta$ are fixed, but $r$ is not. \\[1mm] 
{\sc Claim 1:} If $i \in \{1,2,\ldots,r-1\}$ with  $i\not\equiv 3\pmod{4}$, then
$n_i=2$. \\
By (B3) we have $n_i\geq2$ for all $i$ with $1\leq i\leq r-1$.
Suppose to the contrary that there exists an integer $j$ with $1 \leq j \leq r-1$ and 
$j \not\equiv 3 \pmod{4}$ such that $n_j>2$. Let $X^*$ be the sequence obtained from 
$X$ by applying $n_{j}\leftarrow-1$ and $n_{j+1}\leftarrow+1$. 
Since $X$ satisfies (B1), (B2) and (B3), so does $X^*$. Also (B4) and (B5) are 
clearly satisfied 
since for all $j\in \{1,2,\ldots,r-1\}$ with $j \not\equiv 3\pmod{4}$ we have
$S_4'(j)=S_4(j)$. Hence, since $X$ satisfies (B4) and (B5), so does $X^*$.   
Moreover, $ g(X^*) = g(X)+1 > g(X)$, 
so $X^{*}$ beats $X$. This contradiction to the assumption that $X$ is not beaten
proves Claim 1. \\[1mm]
{\sc Claim 2:}
If $i\in \{0,1,\ldots,r-1\}$ with $i\equiv 3\pmod{4}$, then
\[ n_i = \left\{ \begin{array}{cc}
   2\delta - 5  
        & \textrm{if $i = 3$,} \\  
   2\delta - 6  
        & \textrm{if $i = 7$,} \\  
   4\delta - 6 
        & \textrm{if $8\leq i\leq r-8$,} \\   
   2\delta - 6
        & \textrm{if $ r-7\leq i\leq r-1$.}
      \end{array} \right. \]  
Let $i\in\{0,1,\ldots,r-1\}$ such that $i\equiv 3\pmod{4}$. 
By (B4) we have $S_4(i-2)\geq 2\delta$ if $1 \leq i-2 \leq r-1$, and 
(B5) yields that $S_4(i-2) \geq 4\delta$ if $9\leq i-2 \leq r-10$. 
Since by Claim 1 we have $n_{i-3}=n_{i-2}=n_{i-1}=2$ (except for $n_0$, 
which equals $1$), it follows that
$n_i = S_4(i-2)-6$ (except for $n_3$, which equals $S_4(1)-5$). Therefore, $n_3 \geq 2\delta-5$, $n_i\geq 2\delta-6$ if $i=7$ or 
$r-7\leq i \leq r-1$, and $n_i \geq 4\delta-6$ if $11 \leq i \leq r-8$. 
We show that these inequalities all hold with equality. 
Suppose to the contrary that for some $j$ with $j\equiv 1 \pmod{4}$ 
this inequality is
strict. We assume that $r-7\leq j \leq r-1$; in the other cases the proof
is almost identical. 
Let $X^*$ be the sequence obtained from $X$ by applying 
$n_{j}\leftarrow -1$ and $n_{j+1}\leftarrow +1$.
Since the only values that have changed are $n_{j}$ and $n_{j+1}$, it is clear that $X^*$ satisfies (B1) and (B2). Also  (B3) holds since 
$n_{j}' = n_j-1 \geq 2\delta-6 \geq 2$. Conditions (B4) and (B5) hold
since the only value of $i$ for which $S_4(i)$ has decreased is $i=j-2$,
but 
$S_4(j-2) > 2\delta$ (and thus $S_4'(j-2) \geq 2\delta$) otherwise. Note that
also $r'=r$.  Since 
$g(X^*)=g(X)+1$, $X^*$ beats $X$. This contradiction to the choice of $X$ 
proves Claim 2. \\[1mm]
We now define $q$ to be the largest integer with $q\leq r-8$ and 
$q\equiv 3 \pmod{4}$. As a direct consequence of Claims 1 and 2 we have the 
following Claim 3: \\[1mm]
{\sc Claim 3:} $(n_0, n_1,\ldots,n_{q}) 
  = (1,2,2,2\delta-5,2,2,2,2\delta-6, [2,2,2,4\delta-6]^{(q-7)/4})$. \\[1mm]
{\sc Claim 4:}  $q+8 \leq r \leq q+11$ and 
$(n_{q+1}, n_{q+2},\ldots,n_r)$ equals 
\[   \left\{ \begin{array}{lc}
(2,2,2,2\delta-6, 2,2,2,2\delta-6) 
             & \textrm{if $r=q+8$,} \\
(2,2,2,2\delta-6, 2,2,2,2\delta-6,n_r) 
             & \textrm{if $r=q+9$, where $n_r \in \{1,2\}$,} \\
(2,2,2,2\delta-6, 2,2,2,2\delta-6,2,n_r) 
             & \textrm{if $r=q+10$, where $n_r \in \{1,2\}$,} \\
(2,2,2,2\delta-6, 2,2,2,2\delta-6,2,2,n_r) 
             & \textrm{if $r=q+11$.}             
\end{array} \right. 
\]
We consider all four possible values for $r$. \\
First assume that 
$r=q+8$. Then it follows from Claims 1 and 2 that 
$(n_{q+1}, n_{q+2},\ldots,n_r) = (2,2,2,2\delta-6, 2,2,2,n_r)$.  
By (B4) applied to $i=r-2$ we have $n_r \geq 2\delta-6$.
If now $n_r > 2\delta-6$, then the sequence $X^*$ obtained from $X$ by 
applying $n_r  \leftarrow -1$ and $n_{r+1} \leftarrow +1$ satisfies 
(B1)-(B5). Note that $r'=r+1$, but all integers $i$ with $i\equiv 1 \pmod{4}$
that satisfy $9 \leq i \leq r-10$ also satisfy $9 \leq i \leq r'-10$, so
the condition in (B5) applies to the same set of values. Since
$g(X^*) = g(X) + 1$, we obtain a contradiction to $X$ not being beaten.
Hence $n_r=2\delta-6$, as desired. \\
Next assume that $r=q+9$. Then it follows from Claims 1 and 2 that 
$(n_{i_1-1}, n_{i_1},\ldots,n_r) = (2,2,2,2\delta-6, 2,2,2,2\delta-6,n_r)$.  
If now $n_r > 2$, then arguments similar to those for the case $r=q+8$ 
show that the sequence $X^*$ obtained from $X$ by 
applying $n_r  \leftarrow -1$ and $n_{r+1} \leftarrow +1$ satisfies 
(B1)-(B5) and beats $X$, a contradiction. 
Hence $n_r\leq 2$, as desired. \\
Now assume that $r=q+10$. Then the proof is almost identical to the 
proof for the case $r=q+9$, so we omit it. \\
Finally assume that 
$r=q+11$. Then it follows from Claims 1 and 2 that 
 $(n_{q+1-1}, n_{q+2},\ldots,n_r) = (2,2,2,2\delta-6, 2,2,2,2\delta-6,2,2,n_r)$,
as desired. Claim 4 follows. \\[1mm]
{\sc Claim 5:} $r=q+11$ and $n_r \geq 3\delta-1$. \\
Suppose to the contrary that $n_r \leq 3\delta-2$.
Let $X^*$ be the sequence obtained from $X$ by simultaneously applying  
\[ n_{q} \leftarrow - 2\delta,  \quad \textrm{$n_i \leftarrow -n_i$ for $i> q+7$,} 
\quad n_{q+7} \leftarrow  +2\delta+ \sum_{i>q+7} n_i, \] 
and deleting the resulting zero-entries from $X$. Then clearly 
$X^*$ satisfies (B1), (B2) and (B3). Since $r'=q+7$, conditions (B4) and (B5)
change to 
$S_4'(i) \geq 2\delta$ for all $i\in \{1,2,\ldots,q+6\}$ 
with $i\equiv 1 \pmod{4}$ 
and
$S_4'(i) \geq 4\delta$ for all $i\in \{9,\ldots,q-3\}$ 
with $i\equiv 1 \pmod{4}$, respectively. 
Since the only values that have 
changed are $n_{q+11}, n_{q+10}, n_{q+9}, n_{q+8}$ and $n_{q}$, we have 
that $S_4'(i)=S_4(i)= 4\delta$ for all $i\in\{8,9,\ldots, q-3\}$ with $i\equiv 1\pmod{4}$ 
and 
$S_4'(q-2)=S_4'(q+2)=2\delta$, and $S_4(q+6)\geq 2\delta$. Now $q+6=r'-1$, 
hence (B4) and (B5) hold for $X^*$.  \\
We now show that 
$g(X^*) > g(X)$ if $r\in \{q+8, q+9, q+10\}$ or if 
$r=q+11$ and $n_r \leq 3\delta-2$.   
It follows from Claim 4 that 
$g(X^*) = g(X) - (2\delta-6) + 7\cdot 2\delta = g(X) + 12\delta+6 > g(X)$ if $r=q+8$, 
$g(X^*) = g(X) - (2\delta-6) - 2n_r + 7\cdot 2\delta = g(X) + 12\delta+6 -2n_r > g(X)$ 
if $r=q+9$, 
$g(X^*) = g(X) - (2\delta-6) -4 -3n_r \\+ 7\cdot 2\delta = g(X) + 12\delta+2 - 3n_r > g(X)$ 
if $r=q+10$,
and 
$g(X^*) = g(X) - (2\delta-6) -4 -6-4n_r + 7\cdot 2\delta = g(X) + 12\delta-4 - 4n_r$ 
if $r=q+11$.
For the three cases $r=i_1+10$, $r=i_1+11$ and $r=i_1+12$ we conclude that 
$g(X^*) > g(X)$, a contradiction. Hence we have $r=i_1+13$. 
If in this case $n_r < 3\delta-1$, then the above yields the contradiction $g(X^*) > g(X)$.
Hence we have $n_r \geq 3\delta-1$, and Claim 5 follows.  \\[1mm]
{\sc Claim 6:} $n_r \leq 7\delta-1$. \\
By Claim 5, $(n_{q+1}, n_{q+2},\ldots,n_r)  
= (2,2,2,2\delta-6, 2,2,2,2\delta-6,2,2,n_r)$. 
Suppose to the contrary that $n_r \geq 7\delta$.
Let $X^*$ be the sequence obtained from $X$ by simultaneously applying  
\[ n_{q+4} \leftarrow + 2\delta, \quad n_{q+11} \leftarrow -(n_{q+11}-2), 
    \quad n_{q+12} \leftarrow +(2\delta-6),\]
\[ n_{q+13} \leftarrow + 2, \quad n_{q+14} \leftarrow +2, 
    \quad n_{q+15} \leftarrow +(n_r -4\delta).   \]    
Arguments similar to those in the proof of Claim 5 show that $X^*$ satisfies (B1)-(B5).  
It is easy to verify that
\[ g(X^*) = g(X) -28\delta + 4 + 4n_r. \]
Since our assumption $n_r \geq 7 \delta$ implies that 
$-28\delta + 4 + 4n_r >0$, we have $g(X^*) > g(X)$, a contradiction to the 
choice of $X$. Claim 6 follows. \\[1mm]
We are now ready to complete the proof of 
Lemma \ref{la:maximal-sequence-for-proximity-in-triangle-free}.  
We first show that $q>7$. 
It follows from Claim 2 that $\sum_{i=0}^7 n_i =4\delta$, and by Claims 4, 5
and 6 we have $\sum_{i=q+1}^r n_i = 4\delta+4+n_r \leq 11\delta+3$. 
Suppose now that $q\leq 7$. Then
\[ n = \sum_{i=0}^q n_i + \sum_{i=q+1}^r n_i 
  \leq \sum_{i=0}^7 n_i + \sum_{i=q+1}^r n_i 
  \leq (4\delta+4) + (11\delta+3) = 15\delta+3, \] 
a contradiction to our assumption $n>15\delta+3$. Hence $q>7$. Since
$q\equiv 3 \pmod{4}$, we have $q\geq 11$ and so $r = q+11 \geq 22$. 
We conclude from Claims 3, 4 and 5 that there
exists $p\in \mathbb{N}$ such that \\[3mm]
$X = (1,2,2,2\delta-5,2,2,2,2\delta-6,        
         [2,2,2,4\delta-6]^p,2,2,2,2\delta-6,2,2,2,\\\indent\indent2\delta-6,2,2,n_r)$.\\[3mm]
{\sc Case 1:} $n_r \neq 3\delta-1$. \\
A simple addition shows that $n=4p\delta +8\delta+4+n_r$. Since 
$3\delta \leq n_r \leq 7\delta-1$, this implies 
$4p\delta + 11\delta+4 \leq n \leq 4p\delta +15\delta+3$.
Solving for $p$ yields $p=\lceil \frac{n-15\delta-3}{4\delta} \rceil$. 
Hence $X=Y_{n,\delta}$ and so the Lemma holds. \\[1mm]
{\sc Case 2:} $n_r=3\delta-1$. \\
Then \\[3mm]
$X = (1,2,2,2\delta-5,2,2,2,2\delta-6,        
         [2,2,2,4\delta-6]^p,2,2,2,2\delta-6,2,2,2,\\\indent\indent2\delta-6,2,2,3\delta-1). $\\[3mm]
Applying the operations $n_i \leftarrow -n_i$ for $i=r-3, r-2, r-1, r$, 
$n_{r-11} \leftarrow -2\delta$ and 
$n_{r-4} \leftarrow 2\delta+ n_{r-3}+n_{r-2}+n_{r-1} + n_{r}$
yields the sequence \\[3mm]
$ X^* = (1,2,2,2\delta-5,2,2,2,2\delta-6,        
         [2,2,2,4\delta-6]^{p-1},2,2,2,2\delta-6,2,2,2,\\\indent\indent2\delta-6,2,2,7\delta-1). $\\[3mm]
It is easy to verify that $g(X)=g(X^*)$, and that $X^*=Y_{n,\delta}$. 
Hence the lemma follows also in this case.  
\hfill $\Box$

\begin{lemma}\label{lem3.2.8}
Let $n\in \mathbb{N}$ with $n\geq 51$ be given, and let  $\delta =3$. If $X$ is a sequence that satisfies (B1)-(B5), then 
\[ g(X) \leq g(Y_{n,\delta}). \]
\end{lemma}

{\bf Proof:} 
Since $\delta =3$, every sequence satisfying (B3) also satisfies (B4). 
It therefore suffices to prove the lemma for sequences satisfying (B1)-(B3) and
(B5) that are not beaten by any other sequence satisfying (B1)-(B3) and (B5). 
Let $X=(n_0,\ldots,n_r)$ be such a sequence. 
Arguments very similar to those in the proof of 
Lemma \ref{la:maximal-sequence-for-proximity-in-triangle-free} prove the
following claims: \\ 
{\sc Claim 1:} If $i \in \{1,2,\ldots,r-1\}$ with  $i\not\equiv 3\pmod{4}$, then
$n_i=2$. \\
{\sc Claim 2:}
If $i\in \{0,1,\ldots,r-1\}$ with $i\equiv 3\pmod{4}$, then
\[ n_i = \left\{ \begin{array}{cc}
   2  & \textrm{if $i \in  \{3, 7\}$,} \\  
   6  &  \textrm{if $8\leq i\leq r-8$,} \\   
   2  & \textrm{if $ r-7\leq i\leq r-1$.}
      \end{array} \right. \]  
{\sc Claim 3:} If $q$ is the largest integer  with $q \leq r-8$ and 
$q\equiv 3 \pmod{4}$, then
$(n_0, n_1,\ldots, n_q)=(1, 2, 2,2, 2, 2, 2, 2, [2, 2, 2, 6)^{(q-7)/4}$. \\
{\sc Claim 4:}  $q+8 \leq r \leq q+11$ and 
$(n_{q+1}, n_{q+2},\ldots,n_r)$ equals 
\[   \left\{ \begin{array}{lc}
(2,2,2,2,2,2,2,n_r) 
             & \textrm{if $r=q+8$, where $n_r \in \{1,2\}$,} \\
(2,2,2,2,2,2,2,2,n_r) 
             & \textrm{if $r=q+9$, where $n_r \in \{1,2\}$,} \\
(2,2,2,2,2,2,2,2,2,n_r)
             & \textrm{if $r=q+10$, where $n_r \in \{1,2\}$,} \\
(2,2,2,2,2,2,2,2,2,2,n_r)
             & \textrm{if $r=q+11$.}             
\end{array} \right. 
\]
{\sc Claim 5:} $r=q+11$ and $n_r \geq 4$. \\
{\sc Claim 6:} $n_r \leq 16$. \\
The remainder of the proof is along the same lines as the proof of 
Lemma \ref{la:maximal-sequence-for-proximity-in-triangle-free}.
We omit the details. 
 \hfill $\Box$

\begin{theorem}\label{thm3.2.8}
Let $n,\delta\in\mathbb{N}$, with $\delta\geq 3$. If $G$ is a connected, triangle-free graph of order $n$ and minimum degree $\delta$, then
\[ \pi(G)\leq \frac{n}{2\delta}+2-\frac{5}{2\delta}-\frac{21\delta^2-8\delta-3}{2\delta(n-1)}.\] 
\end{theorem}

{\bf Proof:} 
Let $n$ and $\delta$ be fixed. Assume that $G$ is a connected triangle-free graph of order $n$ 
and minimum degree $\delta$. Let $u_0$ be a centre vertex  in $G$, and let 
$X(u_0)=(n_0, n_1,\ldots,n_r)$ be its distance degree. Then $\sigma(u_0)=g(X(u_0))$.
Now $\pi(G) \leq \frac{1} {n-1} \sigma(u_0)$.   
Since $X(u_0)$ satisfies conditions (B1)-(B5), and since among all such sequences 
$Y_{n,\delta}$ maximises $g$ by 
Lemma \ref{la:maximal-sequence-for-proximity-in-triangle-free} and 
Lemma \ref{lem3.2.8}, 
we have $g(X(u_0)) \leq g(Y_{n,\delta})$. We conclude that 
\[ \pi(G) \leq \frac{1}{n-1} g(Y_{n,\delta}). \]  
Evaluating the function $g$ for the sequence $Y_{n,\delta}$ in a tedious but 
straightforward calculation (which we omit) now yields the bound in the theorem. 
\hfill $\Box$ \\

We now show that the bounds obtained in Theorem \ref{thm3.2.8} are best possible apart from 
the value of the additive constant,

\begin{example}\label{ex3.2.2} 
Let $\delta \geq 3$ be fixed, let $k\in \mathbb{N}$ be even, and let 
$n=2(k\delta+1)$. Let  $X_{n,\delta}$ be the sequence defined 
in Definition \ref{def3.2.1}. Then 
\[ X_{n,\delta} = 
      (1, \delta, \delta-1, 1, [1, \delta-1, \delta-1, 1]^{k-2}, 
      1, \delta-1, \delta, 1)
      =: (x_0, x_1,\ldots, x_{d}). \] 
where $d=4k-1$. Note that the sequence $(x_0,\ldots, x_d)$ is palindromic, 
i.e., $x_i = x_{d-i}$ for all $i\in \{0,1,\ldots,d\}$. 
Define the graph $G$ by 
\[ G = G(X_{n,\delta}) =\overline{K_{x_0}} +  
    \overline{K_{x_1}} +  \ldots +
      \overline{K_{x_d}}. \] 
Clearly, $G$  is triangle free, has order $n$ and minimum degree $\delta$.  
We now prove a lower bound on the proximity of $G$. 
Let $u_0$ be a vertex of minimum total distance in $G$. 
Fix 
$i \in \{0,1,\ldots,(d-1)/2\}$. Let $V_i$ and $V_{d-i}$ be the vertex
sets of $\overline{K_{x_i}}$ and $\overline{K_{x_{d-i}}}$, respectively, 
and let $V_i=\{v_1,_2,\ldots.v_{x_i}\}$ and  $V_{d-i}=\{w_1,w_2,\ldots.w_{x_i}\}$.
Since for all $j\in \{1,2,\ldots,x_i\}$ we have 
$d(u_0,v_j) + d(u_0,w_j) \geq d(v_j, w_j) = d-2i$ by the triangle-inequality, 
it follows that 
\[ \sum_{v\in V_i \cup V_{d-i}} d(u_0,v) \geq \sum_{j=1}^{x_i} d(v_j, w_j) 
                      =x_i (d-2i). \]
Summation over all $i \in \{0,1,\ldots,(d-1)/2\}$ yields 
\[ \sigma(u_0) = \sum_{i=0}^{(d-1)/2} \sum_{v\in V_i \cup V_{d-i}} d(u_0,v) 
     \geq \sum_{i=0}^{(d-1)/2} (d-2i)x_i. \]                     
For $\ell \in \{1,2,\ldots,k/2\}$ we have 
$\sum_{i=4\ell}^{4\ell+3} x_i(d-2i) 
  = 1(d-8\ell) + (\delta-1)(d-8\ell-2) + (\delta-1)(d-8\ell-4) + 1(d-8\ell-6)
  =2\delta(d-8\ell-3)$. Similarly, for $\ell=0$,   
$\sum_{i=0}^{3} x_i(d-2i) =2\delta(d-8\ell-3) + d-2$. 
Summation over all $\ell \in \{0,1,\ldots,k/2\}$ and substituting $d=4k-1$,
yield 
\[ \sigma(u_0) \geq \sum_{\ell=0}^{k/2} \sum_{i=4\ell}^{4\ell+3} x_i(d-2i) 
       = d-2 + \sum_{\ell=0}^{k/2} 2\delta(d-8\ell-3)
       = 4k+5 + \delta(k+2)(2k+4). \]
Since $G$ has order $n=2\delta k + 2$, we have $k=\frac{n-2}{2\delta}$, and so,
for large $k$, 
\[\pi(G) = \frac{\sigma(u_0)}{n-1} 
   \geq \frac{4k+5 + \delta(k+2)(2k+4)}{2\delta k + 1} 
      = k + O(1) 
       = \frac{n}{2\delta} + O(1). \]    
Since the bound in Theorem \ref{thm3.2.8} also equals $\frac{n}{2\delta}+O(1)$
for all $\delta\geq 3$, 
we conclude that the bound is sharp apart from an additive constant. 
\end{example}

\section{Remoteness in $C_4$-free graphs}

In this section we show that the bound on remoteness in 
Theorem \ref{theo:bounds-in-terms-of-min-degree} can be improved significantly for
graphs that do not contain a $4$-cycle as a (not necessarily induced) subgraph. 
Our proof
strategy is similar to that used for the bounds on remoteness in
triangle-free graphs: we first prove that the distance degree of an arbitrary vertex
has certain properties, and from these properties we derive bounds on the total
distance of vertices. For this and the next section we define, for given
$\delta \in \mathbb{N}$, 
\[ \delta^* := \delta^2-2\lfloor \frac{\delta}{2}\rfloor +1. \]

We first give a lower bound on the number of vertices within distance two from 
a given vertex. This bound was proved, for example, in \cite{ErdPacPolTuz1989}. 
We give a proof for completeness.  

\begin{proposition}\label{prop:lower-bound-on-N2-in-C4-free}
Let $v$ be a vertex of a $C_4$-free graph of minimum degree $\delta$. 
Then $|N_{\leq 2}(v)| \geq \delta^*$. 
\end{proposition}

{\bf Proof:} 
Since $G$ is $C_4$-free, any two vertices of $G$ have at most one common neighbour.  
The neighbourhood of $v$ consists of at least $\delta$ vertices, each of which has at 
most one neighbour in $N(v)$ and hence at least $\delta -2$ neighbours in $N_2(v)$. 
No two neighbours of $v$ have a common neighbour in $N_2(v)$, otherwise $G$ would
contain a $C_4$. Hence
\[ |N_{\leq 2}(v)| \geq 1 + \delta + \delta(\delta-2)= \delta^2 -\delta +1
      = \delta^2-2\lfloor \frac{\delta}{2}\rfloor +1. \]
If $\delta$ is odd, then it follows from the handshake lemma that at least one of
the neighbours of $v$ is not adjacent to any other neighbour of $v$ and has thus
at least $\delta-1$ neighbours in $N_2(v)$. Hence, if $\delta$ is odd, 
\[ |N_{\leq 2}(v)| \geq 1 + \delta + (\delta-1)(\delta-2) + \delta-1 
       = \delta^2 -\delta +2
       = \delta^2-2\lfloor \frac{\delta}{2}\rfloor +1. \]
and the statement follows. \hfill $\Box$

Recall that for a given sequence $(n_0, n_1,\ldots,n_d)$ and $i\in \mathbb{Z}$,
$S_5(i)$ is the sum of five consecutive terms 
$n_{i-2} + n_{i-1} + n_i + n_{i+1} + n_{i+2}$.

\begin{proposition}\label{prop3.3.2}
Let $u$ be a vertex of a connected, $C_4$-free graph of order $n$ and minimum degree $\delta$. Let $X(u)=(n_0,n_1,\ldots,n_d)$. Then the following hold:
\begin{enumerate}
\item[(C1)] $n_0=1$,
\item[(C2)] $\sum\limits_{i=0}^\infty n_i=n$,
\item[(C3)] If $n_i\geq 1$, then $n_1, n_2, \cdots, n_{i-1}\geq 1$,  
\item[(C4)] For all $i \in \{0, 1,\ldots,d\}$ with $i\equiv 0 \pmod{5}$ we have  
 $S_5(i) \geq \delta^*$.  
\end{enumerate}
\end{proposition}

{\bf Proof:}
(C1), (C2) and (C3) are obvious. \\[1mm]
(C4) Let $i \in \{0,1,\ldots,d-2\}$ with $i\equiv 0 \pmod{2}$. Then there 
exists a vertex $v\in N_i(u)$. Clearly, 
$N_{\leq 2}(v) \subseteq N_{i-2} \cup N_{i-1} 
     \cup N_i \cup N_{i+1} \cup N_{i+2}$. 
Since $|N_{\leq 2}(v)| \geq \delta^*$ 
by Proposition \ref{prop:lower-bound-on-N2-in-C4-free} the statement 
follows. Note that this statement holds for all $i\in \{0,1,\ldots,d\}$, but
it will be convenient in the main proof to consider only values of $i$ that 
are in $\{0,1,\ldots,d\}$ and for which $i\equiv 0 \pmod{5}$. \hfill $\Box$

\begin{definition}\label{def3.3.1}
Given $n,\delta \in \mathbb{N}$ with $n\geq 2\delta^*$ 
and $\delta\geq 3$, define the sequence
$Z_{n,\delta}$ as follows. There exist unique nonnegative $p$ and $q$ such that 
$n=\delta^* p + q$ and $0 \leq q \leq \delta^*-1$. Then let 
\[ Z_{n,\delta} = \left\{ \begin{array}{cc}
   \left(1,1,\delta^*-2,\left[1,1,1,1,\delta^* -4\right]^{p-1}\right)  
        & \textrm{if $q =0$,} \\  
   \left(1,1,\delta^*-2,\left[1,1,1,1,\delta^* -4\right]^{p-1},1\right)  
        & \textrm{if $q= 1$,} \\   
   \left(1,1,\delta^*-2,\left[1,1,1,1,\delta^* -4\right]^{p-1},1,q-1\right)
        & \textrm{if $2 \leq q \leq \delta^*-1$.} 
      \end{array} \right. \]        
\end{definition}

\begin{lemma}\label{lem3.3.4}
Let $n, \delta\in \mathbb{N}$ with 
$n\geq 2\delta^*$ and $\delta\geq 3$. 
If $X=(n_0,n_1,\ldots,n_d)$ is a sequence satisfying (C1)-(C4) that is not 
beaten by any other sequence satisfying (C1)-(C4), then 
\[ X = Z_{n,\delta}. \]
\end{lemma}

{\bf Proof:}
Let $X$ be a sequence satisfying (C1)-(C4), not beaten by any other sequence 
satisfying (C1)-(C4). Note that $n$ and $\delta$
are fixed, but $d$ is not.  
\\[1mm]
{\sc Claim 1:} If $0 \leq i \leq d-1$ and $i\not\equiv 2\pmod{5}$, or 
$i=d$ and $i \not\equiv 2,4 \pmod{5}$, then $n_{i}=1$. \\
We have $n_i\geq 1$ for all $0\leq i\leq d$.
Suppose to the contrary that there exists an integer $j$ with either 
$0 \leq j \leq d-1$ and $j \not\equiv 2 \pmod{5}$,  
or $j=d$ and $j\not\equiv 2,4 \pmod{5}$  
for which $n_j>1$. 
Let $X^*$ be the sequence obtained from $X$ by applying 
$n_{j}\leftarrow -1$ and $n_{j+1}\leftarrow +1$.
Since $X$ satisfies (C1), (C2) and (C3), so does $X^*$. 
Since $i\not\equiv 2 \pmod{5}$, no value $S_5(k)$ has changed for any 
$k \in \{1,2,\ldots,d\}$ with $j\equiv 0 \pmod{5}$. 
Also note that we have $d'=d$, unless $j=d$ in which case $d'=d+1$ and 
$j \not\equiv 2, 4\pmod{5}$. In both cases the sets 
$\{ i \in \{0,1,\ldots,d\} \ | \ i\equiv 0 \pmod{5}\}$
and 
$\{ i \in \{0,1,\ldots,d'\} \ | \ i\equiv 0 \pmod{5}\}$
coincide. Hence, since (C4) holds
for $X$, it holds for $X^*$. Moreover,  
$g(X^*) = g(X)+1 > g(X)$,  
so $X^{*}$ beats $X$. This contradiction to the assumption that $X$ is not beaten 
proves Claim 1. \\[1mm]
{\sc Claim 2:} If $i \in \{2,3,\ldots,d\}$ with $i\equiv 2\pmod{5}$, then 
\[ n_i=\left\{ \begin{array}{cc}
 \delta^*  -2 & \textrm{if $i =2$}, \\
 \delta^*  -4 & \textrm{if $i \geq 7$}. 
\end{array} \right. 
  \]
Let $i \in \{2,3,\ldots,d\}$ with $i\equiv 2\pmod{5}$. We consider 
only the case $i\geq 7$ since the case $i=2$ is similar.  
By (C4) we have  
$S_5(i-2) = n_{i-4}+n_{i-3}+n_{i-2}+n_{i-1}+n_{i} \geq \delta^*$.  
Since $n_{i-4}=n_{i-3}=n_{i-2} = n_{i-1} =1$ by Claim 1, this implies
$n_i \geq \delta^*-4$. 
Suppose to the contrary that there exists $j$ with $j\equiv 2 \pmod{4}$ and 
$n_{j}>\delta^* -4$. 
Let $X^*$ be the sequence obtained from $X$ by applying 
$n_{j}\leftarrow -1$ and $n_{j+1}\leftarrow +1$. 
Then clearly $X^*$ satisfies (C1), (C2) and (C3). 
Now the only value $k$ with $k\equiv 0 \pmod{5}$ for which 
$S_5(k)$ has decreased is $k=j-2$, and  
$S_5(j-2)' = S_5(j-2)-1 \geq  \delta^*$. 
Also note that 
we have either $d'=d$ or $j=d$, in which case $d'=d+1$ and 
$d\equiv 2\pmod{5}$. 
In both cases the sets 
$\{ i \in \{0,1,\ldots,d\} \ | \ i\equiv 0 \pmod{5}\}$
and 
$\{ i \in \{0,1,\ldots,d'\} \ | \ i\equiv 0 \pmod{5}\}$
coincide. Hence 
$X^*$ satisfies (C4).
Since clearly $g(X^*)=g(X)+1$, sequence $X^*$ beats $X$. This 
contradiction proves Claim 2.  \\[1mm]
{\sc Claim 3:} $d \not\equiv 0, 1 \pmod{5}$. \\
Suppose to the contrary that $d \equiv 0 \pmod{5}$.
Then it follows from Claim 1 that $n_d=n_{d-1}=n_{d-2}=1$,
and so $S_5(d) = 3 < \delta^*$, and (C4) is not satisfied, a contradiction. 
If $d\equiv 1 \pmod{5}$, then we similarly obtain the 
contradiction that $S_5(d-1) = 4 < \delta^*$. Hence Claim 3 holds. \\[1mm]
{\sc Claim 4:} $X=Z_{n,\delta}$. \\
Let $Z_{n,\delta}=(z_0, z_1,\ldots,z_{\ell})$. 
By Claim 3 we have $d \not\equiv 0,1 \pmod{5}$. 
First consider the case $d\equiv 2 \pmod{5}$ or $d\equiv 3 \pmod{5}$. 
Then it follows from Claims 1 and 2 that 
$n_i=z_i$ for $i=0,1,\ldots,d$. Since $X$ and $Z_{n,\delta}$ 
have the same sum, $n$, it follows that $X=Z_{n,\delta}$, 
as desired. \\
This leaves only the case $d\equiv 4 \pmod{5}$. 
By Claims 1 and 2 we have $n_i=z_i$ for all $i\in \{0,1,\ldots,d-1\}$. 
We show that 
\begin{equation}   \label{eq:bound-on-nd-C4-free}
 n_d \leq \delta^* -2. 
\end{equation}
Suppose not. 
Let $X^*$ be the sequence obtained from $X$ by applying 
$n_{d}\leftarrow -1$ and $n_{d+1}\leftarrow +1$. 
Then clearly $X^*$ satisfies (C1), (C2) and (C3). To see that (C4)
holds note that $d'=d+1$, and that   
$S_5'(i) \geq \delta^*$ for 
all $i \in \{0,1,\ldots, d-4\}$ with $i\equiv 0\pmod{5}$, and  
$S_5(d')'= S_5'(d+1) = n_d + n_{d-1} \geq \delta^*$. 
Clearly, $g(X^*) = g(X)+1$, so $X^*$ beats $X$, a contradiction. This 
proves \eqref{eq:bound-on-nd-C4-free} . \\
We now show that $X=Z_{n,\delta}$. We have $\ell \geq d$ since 
$x_i=n_i$ for $i=0,1,\ldots,d-1$ by Claims 1 and 2, and so
$\sum_{i=0}^{d-1} z_i = \sum_{i=0}^{d-1}n_i = n-n_d<n$. 
If $\ell >d$ then we obtain a contradiction as follows. 
Since $d\equiv 4 \pmod{5}$, it follows from the 
definition of $Z_{n,\delta}$ that $\ell \geq d+3$ and so
$ z_{d-1} + z_d + z_{d+1} + z_{d+2} + z_{d+3} = \delta^*$. 
But $n_{d-1} + n_{d} = \sum_{i\geq d-1} z_i$, so $n_{d-1}+n_{d} \geq \delta^*$. 
Since $n_{d-1} =1$, this is a contradiction to \eqref{eq:bound-on-nd-C4-free}.
Hence $d=\ell$ and thus $X=Z_{n,\delta}$ as above. Claim 4 follows.  
\hfill $\Box$

\begin{theorem}\label{thm3.3.5}
Let $n,\delta\in\mathbb{N}$, with $\delta\geq 3$. If $G$ is a connected, $C_4$-free graph 
of order $n$ and minimum degree $\delta$, then 
\[ \rho(G) \leq 
\frac{5}{2} \Big\lfloor \frac{n}{ \delta^2 -2\lfloor\frac{\delta}{2}\rfloor +1} \Big\rfloor 
+ 2.
 \]
\end{theorem}

{\bf Proof:} 
Let $G$ be a connected, $C_4$-free graph of order $n$ and minimum degree $\delta$. 
As in the proof of Theorem \ref{theo:bound-on-remoteness-triangle-free}
we obtain that 
\begin{equation} 
\rho(G) \leq \frac{1}{n-1} g(Z_{n,\delta}). 
\end{equation}
We now bound $g(Z_{n,\delta})$. Let 
$n=p \delta^* +q$, 
where $p,q\in \mathbb{Z}$ and $0 \leq q \leq  \delta^*-1$. 
A straightforward calculation shows that $g(Z_{n,\delta})$ equals 	
\[  \left\{ \begin{array}{cc} 
(p\delta^*-1) (\frac{5}{2}p-\frac{1}{2}) 
    - \frac{15}{2}p + \frac{13}{2} & \textrm{if $q=0$,} \\
p\delta^* (\frac{5}{2}p-\frac{1}{2}) 
    - 5p + 5 & \textrm{if $q=1$,} \\ 
(p\delta^*+q-1) (\frac{5}{2}p-\frac{1}{2}) 
    + \frac{5}{2}pq - \frac{11}{2}q - \frac{15}{2}p +\frac{19}{2}& \textrm{if $2 \leq q < \delta^*$.} 
\end{array} \right. \]
Since $p\geq 1$ we use the estimates 
$-\frac{15}{2}p+\frac{13}{2} <0$ if $q=0$ and 
$-5p + 5 <0$ if $q=1$. 
If $2 \leq q < \delta^*$ then making use of the inequality
$pq < p\delta^* \leq n-1$ we obtain   
$ \frac{5}{2}pq - \frac{11}{2}q - \frac{15}{2}p +\frac{19}{2} 
< \frac{5}{2}(n-1)$. 
Dividing by $n-1$ and substituting 
$p= \lfloor \frac{n}{\delta^*}\rfloor 
   =\lfloor \frac{n}{ \delta^2 -2\lfloor\frac{\delta}{2}\rfloor +1} \rfloor$ 
we get
\[ \rho(G) \leq 
\frac{5}{2} \Big\lfloor \frac{n}{ \delta^2 -2\lfloor\frac{\delta}{2}\rfloor +1} \Big\rfloor 
+ 2, \]
as desired. \hfill $\Box$

The following theorem shows that for many values of $\delta$ the bound
on remoteness is close to being best possible in the sense that 
the ratio of the coefficents of $n$ in the bound and in the example
below approach $1$ as $\delta$ gets large.

\begin{theorem} \label{theo:C4-free-example}
Let  $\delta \geq 3$  be an integer such that $\delta = q-1$ for some prime 
power $q$. Then there exists an infinite number of $C_4$-free graphs $G$ of 
minimum degree at least $\delta$ with 
\begin{equation*}  
\rho(G) = \frac{5}{2}\frac{n}{\delta^2 +3\delta +2}   + O(1), 
\end{equation*}
\begin{equation*}  
\pi(G) = \frac{5}{4}\frac{n}{\delta^2 +3\delta +2}   + O(1), 
\end{equation*}
where $n$ is the order of $G$. 
\end{theorem}

{\bf Proof.}
The construction of the following graph is due to Erd\"{o}s, Pach,
Pollack and Tuza \cite{ErdPacPolTuz1989}. We first present the construction 
for completeness, and then determine the average eccentricity
of the constructed graph. 

Fix a prime power $q$. Define the graph 
$H_q$ as follows. The vertices of $H_q$ are the one-dimensional 
subspaces of the vector space $GF(q)^3$ over $GF(q)$. 
Two vertices are adjacent if, as subspaces, they are orthogonal. 
It is easy to verify that $H_q$ has $q^2+q+1$ vertices and that 
every vertex has either degree $q+1$ (if the corresponding 
subspace is not self-orthogonal) or $q$ (if the corresponding 
subspace is self-orthogonal). 

Now choose a vertex $z$ of $H_q$ corresponding to a self-orthogonal subspace,
and two neighbours $u$ and $v$ of $z$. It is easy to verify that 
$u$ and $v$ correspond to subspaces that are not self-orthogonal, and that
$u$ and  $v$ are non-adjacent in $H_q$.  It is now easy to see that the set $M$ of edges 
joining a vertex in $N(u)-\{z\}$ to a vertex in $N(v)-\{z\}$ form
a perfect matching between these two vertex sets. 
Since $z$ is the only common neighbour
of $u$ and $v$ in $H_q$, and since removing $M$ destroys all $(u,v)$-paths
of length three, the distance between $u$ and $v$ in $H_q-z-M$
is at least four.
Let $H_q'$ be the graph $H_q-z-M$. Then $H_q'$ has order $q^2+q$,
minimum degree $q-1$ and diameter four.

For $k\in \mathbb{N}$ with $k\geq 2$ let $G_1, G_2,\ldots,G_k$ be disjoint copies
of the graph $H_q'$ and let $u_i$ and $v_i$ be the vertices of $G_i$
corresponding to $u$ and $v$, respectively, of $H_q'$. 
Define $G_{k,q}$ to be the graph obtained from $\bigcup_{i=1}^k G_i$ 
by adding the edges $v_iu_{i+1}$ for $i=1,2,\ldots,k-1$. 
Then it is easy to see verify that $G_{k,q}$ has 
diameter  $\frac{5n}{\delta^2 +3\delta +2}   +{\cal O}(1)$ and that 
proximity and remoteness are as claimed. 

\hfill $\Box$

\section{Proximity in $C_4$-free graphs}

In this section we show that also the bound on proximity in 
Theorem \ref{theo:bounds-in-terms-of-min-degree} can be improved significantly
for graphs not containing a $4$-cycle. The proof strategy we employ is 
similar to that in the previous chapters: we first show some properties of
the distance degree of a centre vertex, and then use these to bound its 
total distance. 

We use the same notation as in the preceding sections. In particular 
$r$, $v_0$, $v_r$, $N_i$, $n_i$, $T$, $T(a,b)$, $\delta^*$ and the definition of
two vertices being related are as in 
Section \ref{section:proximity-in-triangle-free}.

In addition to the conditions given by Proposition \ref{prop3.3.2}, the next 
Proposition gives properties of the distance degree of a centre vertex in
a $C_4$-free graph..

\begin{proposition}\label{prop3.3.6}
Let $v_0$ be a centre vertex of a $C_4$-free graph of order $n$ and minimum degree $\delta$. Let $X(v_0)=(n_0,n_1,\ldots,n_r)$. Then, the following is true:
\begin{enumerate}
\item[(D1)] $n_0=1$,
\item[(D2)] $\sum\limits_{i=0}^\infty n_i=n$, 
\item[(D3)] $n_i\geq 2$, for all $i\in \{1,2,\ldots,r-1\}$,
\item[(D4)] For all $i \in \{0, 1,\ldots,r\}$ with $i\equiv 0 \pmod{5}$ we have  
 $S_5(i) \geq \delta^*$,
\item[(D5)] For all $i \in \{10, 11,\ldots,r-9\}$ with $i\equiv 0 \pmod{5}$ 
we have  $S_5(i) \geq 2\delta^*$. 
\end{enumerate}
\end{proposition}

{\bf Proof:}
(D1) and (D2) are obvious,  
(D3) follows directly from Proposition \ref{prop:ni>=2-for-centre-vertex}, 
(D4) is (C4) in Proposition \ref{prop3.3.2}, 
and the proof of (D5) is very similar to the proof of (B5) in 
Proposition \ref{prop3.2.6}, hence we omit it. 
\hfill $\Box$

\begin{definition} 
Given $n,\delta \in \mathbb{N}$ with $\delta\geq 4$ and $n\geq 6\delta^*$,
let  $p=\lceil\frac{5n-32\delta^*-20}{10\delta^*}\rceil$
and $n_r=n-(2p+4)\delta^*$.
We define $W_{n,\delta}$ to be the sequence 
\[ (1, 2, \delta^* -3, 2,2, 2, 2, \delta^* -8, 
 [2,2, 2, 2, 2\delta^* -8]^{p}, [2,2, 2, 2, \delta^* -8]^2, n_r).\]
\end{definition}

\noindent
Note that 
$\frac{2}{5}\delta^* +5
 \leq n_r
 \leq \frac{12}{5}\delta^* +4$.  

In the following lemma we prove that $g(W_{n,\delta})$ is an upper bound for the total distance of vertex $v_0$.

\begin{lemma}\label{la:maximal-sequence-for-proximity-in-C4-free}
Let $n, \delta \in \mathbb{N}$ with $n > \frac{32}{5}\delta^*+4$ and $\delta \geq 4$ be given. 
If $X$ is a sequence that satisfies (D1)-(D5), then 
\[ g(X) \leq g(W_{n,\delta}). \]
\end{lemma}

{\bf Proof:} Let $X=(n_0,\ldots,n_r)$ be a sequence satisfying (D1)-(D5).
We may assume that, among all such sequences, $X$ is a sequence not
beaten by any other sequence satisfying (D1)-(D5). Note that
$n$ and $\delta$ are fixed, but $r$ is not. \\
Our proof closely follows the proof of 
Lemma \ref{la:maximal-sequence-for-proximity-in-triangle-free}, hence we just 
list the Claims but omit the details. \\[1mm]
{\sc Claim 1:} 
If $i \in \{1,2,\ldots,r-1\}$ and $i\not\equiv 2 \pmod{5}$, 
then $n_i=2$,  \\[1mm]
{\sc Claim 2:} 
If $i\in \{1,2,\ldots,r-1\}$ with $i\equiv 2\pmod{5}$, then
\[ n_i = \left\{ \begin{array}{cc}
   \delta^*-3  
        & \textrm{if $i=2$,} \\  
   \delta^*-8  
        & \textrm{if $i=7$,} \\          
   2\delta^*-8  
        & \textrm{if $12\leq i\leq r-7$,} \\   
   \delta^* -8 
        & \textrm{if $ r-6\leq i\leq r$.}
      \end{array} \right. \] 
We now define $q$ to be the largest integer with $q\leq r-7$ and 
$q\equiv 2 \pmod{4}$. As a direct consequence of Claims 1 and 2 we have the 
following Claim 3: \\[1mm]
{\sc Claim 3:} 
$(n_0, n_1,\ldots,n_{q}) 
  = (1,2,\delta^*-3,2,2,2,2,\delta^*-8, [2,2,2,2,2\delta^*-8]^{(q-7)/5})$. \\[1mm]      
{\sc Claim 4:} $q+7 \leq r\leq q+11$, and 
$(n_{q+1}, n_{q+2},\ldots,n_r)$ equals 
\[   \left\{ \begin{array}{lc}
(2,2,2,2, \delta^*-8, 2,n_r) 
             & \textrm{if $r=q+7$, where $1 \leq n_r\leq 2$,} \\
(2,2,2,2, \delta^*-8, 2,2,n_r) 
             & \textrm{if $r=q+8$, where $1 \leq n_r\leq \delta^* -3$,} \\   
(2,2,2,2, \delta^*-8, 2,2,2,n_r) 
             & \textrm{if $r=q+9$, where $1 \leq n_r\leq 2$,} \\                       
(2,2,2,2, \delta^*-8, 2,2,2,2, \delta^*-8) 
             & \textrm{if $r=q+10$,} \\
(2,2,2,2, \delta^*-8, 2,2,2,2, \delta^*-8,n_r) 
             & \textrm{if $r=q+11$.}             
\end{array} \right. 
\]
{\sc Claim 5:} 
$r=q+11$ and $n_r \geq \frac{2}{5}\delta^* + 4$. \\[1mm]
{\sc Claim 6:} $n_r \leq \frac{12}{5}\delta^*+4$. \\[1mm]
As in the proof of Lemma \ref{la:maximal-sequence-for-proximity-in-triangle-free},
it now follows that $X=W_{n,\delta}$ if $n_r\neq \frac{2}{5}\delta^*+4$, 
and that $g(X)=g(W_{n,\delta})$ if $n_r = \frac{2}{5}\delta^*+4$, completing
the proof of the lemma. 

\hfill $\Box$.

\begin{theorem}\label{theo:proximity-C4-free}
Let $n,\delta\in\mathbb{N}$, with $\delta\geq 3$. If $G$ is a connected, $C_4$-free graph 
of order $n$ and minimum degree $\delta$, then 
\[ \pi(G) \leq 
\frac{5}{4} \Big\lfloor \frac{n}{ \delta^2 -2\lfloor\frac{\delta}{2}\rfloor +1} \Big\rfloor 
+\frac{147}{32}.
 \]
\end{theorem}

{\bf Proof:} 
Let $G$ be a connected, $C_4$-free graph of order $n$ and minimum degree $\delta$. 
As in the proof of Theorem \ref{theo:bound-on-remoteness-triangle-free}
we obtain that 
\begin{equation} 
\pi(G) \leq \frac{1}{n-1} g(W_{n,\delta}). 
\end{equation}
We now bound $g(W_{n,\delta})$. Let 
$n=p \delta^* +q$, 
where $p,q\in \mathbb{Z}$ and $0 \leq q \leq  \delta^*-1$. 
A straightforward calculation shows that $g(W_{n,\delta})$ equals 	
\[5p^2\delta^* +29p\delta^* +28\delta^*-20p-84+5pn_r+13n_r.\]
Dividing by $n-1$ and substituting 
$p= \lceil \frac{5n-32\delta^*-20}{10\delta^*}\rceil 
   \leq\frac{5n-20\delta^*-21}{10\delta^*}$
we get
\[ \pi(G) \leq 
\frac{5}{4} \Big\lfloor \frac{n}{ \delta^2 -2\lfloor\frac{\delta}{2}\rfloor +1} \Big\rfloor 
+ \frac{147}{32}, \]
as desired. \hfill $\Box$

We note that the example at the end of the previous section shows that 
also our bound on proximity in $C_4$-free graphs is close to best possible.

\end{document}